\documentclass[11pt,a4paper,reqno]{amsart}
\usepackage{graphicx}
\graphicspath{ {./images/} }
\usepackage{epsfig}
\usepackage{cite}
\usepackage{color}
\usepackage{float}
 \newtheorem{thm}{Theorem}[section]
 \newtheorem{cor}[thm]{Corollary}

 \theoremstyle{definition}
 
 \theoremstyle{remark}
 \newtheorem{rem}[thm]{Remark}

  \newcommand{\E}{\mathbf{E}}
  \newcommand{\cR}{\mathcal{R}}

 \newcommand{\norm}[1]{\left\Vert#1\right\Vert}
 
\begin{document}
\setcounter{page}{1}
\begin{flushleft}
{\scriptsize Appl. Comput. Math., V.xx, N.xx, 20xx, pp.xx-xx}
\end{flushleft}
\bigskip
\bigskip
\title[S.B. Kologrivova,  E.A. Pchelintsev: Improved parameter estimation ... ] {Improved parameter estimation for a family of exponential distributions}
\author[Appl. Comput. Math., V.xx, N.xx,  20xx]{S.\,B. Kologrivova$^1$,\,  E.\,A. Pchelintsev$^1$
 }
\thanks{$^1$Department of Mathematical Analysis and Theory of Functions, Tomsk State University, Lenina av., 36, 634050, Tomsk, Russia,
\\ \indent\,\,\,e-mail: skologrivova@gmail.com
\\ \indent
  \em \,\,\,Manuscript received xx}

\begin{abstract}
In this paper, we consider the problem of parameter estimating for a family of exponential distributions. We develop the improved estimation method, which generalized the James--Stein approach for a wide class of distributions. The proposed estimator dominates the classical maximum likelihood estimator under the quadratic risk. The estimating procedure is applied to special cases of distributions. The numerical simulations results are given. 

\bigskip
\noindent Keywords: Exponential distributions, dependent observations, maximum likelihood estimator, shrinkage estimator, improved estimation, quadratic risk, numerical simulation

\bigskip \noindent AMS Subject Classification: 62F10, 62E17

\end{abstract}
\maketitle

\smallskip
\section{Introduction} 

Let $X=(X_1,...,X_d)$ be a vector of random variables, the distribution of the $i$-th component of which is given by the density 
\begin{equation}\label{den}
    f_{\theta_i}(x)=\exp\{\theta_i b(x)-\psi(\theta_i)\}k(x) \,, \quad x\in\mathbb{R}\,,\quad i=1,...,d\,,\
\end{equation}
where  $b\,,\psi\,,k$ are $\mathbb{R} \to \mathbb{R}$ functions providing simple regularity conditions, $\theta=(\theta_1,...,\theta_d)$ is a vector of unknown parameters from some bounded set $\Theta \subset \mathbb{R}^d$.

The problem is to estimate the vector $\theta$ from observations $X$ under quadratic risk 
\begin{equation}\label{risk}
    \cR(\theta,\delta)=\E_{\theta}\norm{\theta-\delta}^2\,, \quad \norm{x}^2:=\sum_{i=1}^d x_i^2\,.
\end{equation}
Here $\delta$ is some estimator (measurable function of observations, i.e. $\delta=\delta(X)$).

The exponential distributions play important role in Statistics and its applications. For example, besides numerous applications of the normal distribution, the random variables possessing a distribution~\eqref{den} are used in  generalized linear models for functional imaging \cite{Friston}, for approximation of cell interdivision times \cite{Golubev}, for approximation of cancer incidence \cite{Belikov}, etc.

There are a number of methods for point estimation of unknown parameters, including the method of moments and the maximum likelihood method. The maximum likelihood estimator (MLE) has been considered best-in-class for long time due to consistency and asymptotic efficiency \cite{LehmannCasella1998}. 

Ch. Stein in \cite{7} considering the estimation problem of the mean $\theta$ of a $d$-variate normal distribution with the identity covariance matrix, was shown that in this case MLE $\delta^0=X$ is inadmissible if $d\geq 3$, i.e. there exists some another estimator which outperforms MLE in mean square accuracy sense. A specific type of such estimator was introduced by James and Stein in  \cite{6} :
\begin{equation}
    \delta^{JS}(X)=\left(1-\frac{d-2}{{\|X\|}^2}\right)X
\end{equation}
which is much superior in risk to $X$ in any case in which $\norm{\theta}$ is small.
This result raised the interest of scientists, led to a long discussion and development of the improved estimation theory. 
Many authors have extended the inadmissibility results to cover estimation of the location parameter of a wide class of location invariant distributions for a broad range of loss functions. In a very wide class of problems, then, the best invariant estimator is inadmissible (see, for example, \cite{Berger1975, Brown1975}, and references therein).
Today the James-Stein's estimate is well known among statisticians \cite{LehmannCasella1998}, econometricians \cite{GreenbergWebster1998} and engineering community \cite{MKP1998, XSS1994}.
However, James-Stein's estimate applies only to a normal distribution. The generalization of these results to the case of a broader class of distributions is a relevant subject.
There are many papers in which the authors extended the results of Stein to spherically symmetric distributions (see \cite{BrandweinStrawderman2012} and references therein). Recently, in \cite{Pchelintsev2013, KPP2014} a modified James--Stein type estimators for the parameters in a regression model with conditionally-Gaussian noise were constructed. Such estimates dominate the MLE for $d \geq 2$, when the advantages of the James-Stein estimate are only apparent for $d \geq 3$. Moreover proposed method allows to control the quadratic risk. 
A natural extension of the results of improved estimation in Gaussian observation schemes is to consider the problem of parameters estimating of distributions from an exponential family, which contains, in particular, the normal distribution.
Hudson \cite{Hudson1978} was one of the first to address the issue of improved estimation of parameters in exponential distributions  of the form \eqref{den} from independent observations. 
His approach was based on constructing a James-Stein type estimator with similar quadratic risk. The following estimator was proposed
\begin{equation}\label{Hudson_est}
        \delta^{H}=X-\frac{(d-2)}{S}b(X)\,,
    \end{equation}
where $S= \sum\limits_{i=1}^{d} b^2(X_i)$.
To prove the improvement for the estimator $\delta^{H}$ one needs to estimate the quadratic risks difference 
$\cR(\theta,\delta^{H})-\cR(\theta,\delta^{0})$ frome above by some non-positive quantity. As it was shown by Hudson, this is true if the densities \eqref{den} satisfy the analog of Stein's lemma conditions. For this it's necessary that in \eqref{den} the function $b(x)=\int (a(x))^{-1} dx$ is the indefinite integral, where $a$ is a positive function and $\int (a(x))^{-1} dx$ exists in the interior of the domain $X_i$. Then
 the following identity holds
\begin{equation}\label{lemStein}
    \E((X_i-\theta_i)g(X_i))=\E(a(X_i)g'(X_i))\,
\end{equation}
for any absolutely continuous real valued function $g$ such that $\E |a(X_i)g'(X_i)|<\infty$.
 Hence, the estimate \eqref{Hudson_est} dominates the MLE when $d\geq 3$.

It should be noted that subclass of the continuous exponential family \eqref{den}--\eqref{lemStein} contains in particular the normal ($\theta$, 1), gamma ($\theta$, 1) and scaled chi-squared distributions.

 Later Berger in \cite{2}  and Ghosh et al. in \cite{3} constructed another improved estimators in the case of exponential distributions. But in this papers the authors did not receive explicit expressions for risks. This makes it difficult to study the properties of proposed estimates in practice.

The goal of this paper is to construct a new improved estimator for the parameters of exponential distributions mixing the approaches from \cite{Hudson1978} and \cite{Pchelintsev2013}. As it will be shown below such method has some advantages. In particular, it allows to consider the estimation problem for dependent observations schemes.

In Sect. 2, we propose the estimates which dominate MLE under quadratic risk. The main result is done in Theorem~\ref{Th1}. Sect. 3 presents examples of estimates for some special cases of distributions. Here we also apply the improved estimation method for classical exponential distribution for which the identity \eqref{lemStein}  fails.
Sect. 4 presents the results of numerical simulations for comparing the empirical risks of three estimates: MLE, the James-Stein / Hudson estimate, and the proposed estimates  in Sect. 3.

\bigskip
\section{Main results}

Now we define the shrinkage estimate for the unknown vector $\theta$ in the following form
\begin{equation}\label{impr}
    \delta^{*}=X - \frac{c}{\sqrt{S}}b(X)\,,
\end{equation}
where $c$ is some positive constant depending of $d$, statistics $S=\sum_{i=1}^d b^2(X_i)$ and $b(X)=(b(X_1),...,b(X_d))$.

To study the properties of estimates \eqref{impr} one needs the following additional condition\\
$\mathbf{C}$: {\it there exists a positive $a^*$ such that $\E_{\theta}(S^{-1/2})\geq a^*$.}\\
This condition means that the second term in \eqref{impr} is separated from zero in mean, and hence the estimate $\delta^*$ is different of MLE.

\begin{thm}\label{Th1}
    Let the components of the vector $X$ have densities \eqref{den}--\eqref{lemStein} and the condition $\mathbf{C}$ holds. Then for all $d \geq 2$ \\
   a) $\cR(\theta,\delta^*)=\cR(\theta,\delta^0)+c^2-2c(d-1)\E_{\theta}(S^{-1/2})$;\\
   b) the estimate \eqref{impr} dominates the MLE and is minimax for all $0<c<2(d-1)a^{*}$;\\
   c) the uniformly optimal choice of $c$ is 
   $c=(d-1)a^*$.
\end{thm}

{\bf Proof.}
	The risk of MLE $\delta^0 = X$: 
 \begin{equation}
      \cR(\theta,\delta^0)=\E_{\theta} \sum\limits_{i=1}^d (X_i-\theta_i)^2\,.
 \end{equation}
Now we represent the $i$-th component of estimate \eqref{impr} as $\delta^{*}_i= X_i + g_i(X)$, where $g_i(x):=-c b(x_i)/\sqrt{S}$. Then its risk 
\begin{gather*}
    \cR(\theta,\delta^{*})=\E_{\theta}\sum\limits_{i=1}^d(X_i+g(X_i)-\theta_i)^2= \\
    =\E_{\theta}\sum\limits_{i=1}^d (X_i-\theta_i)^2+2\E_{\theta}\sum\limits_{i=1}^d (X_i - \theta_i)g_i(X) +\E_{\theta}\sum\limits_{i=1}^d g^2_i(X).
\end{gather*}
One notes that the first term is the risk of MLE and the last term is equal to $c^2$. Then for the difference of the quadratic risks of estimate \eqref{impr} and MLE we have
$$
\Delta(\theta)=\cR(\theta,\delta^*)-\cR(\theta,\delta^{0})=c^2+2\sum\limits_{i=1}^d \E_{\theta}(X_i - \theta_i)g_i(X)\,.
$$ 
Using here the identity \eqref{lemStein}, we get
\begin{equation}\label{riskmain}
		\Delta (\theta) = c^2+2\sum\limits_{i=1}^d \E_{\theta} [a(X_i)g'_{i}(X)]\,.
\end{equation}
It's easy to check that the derivative
\begin{equation}\label{deriv_g}
		g'_{i}(x)\bigl|_{x=X_i}=\frac{\partial g_i(x)}{\partial x_i}\bigl|_{x=X_i}=-\frac{c b'(X_i)}{\sqrt{S}}+\frac{c b^2(X_i) b'(X_i)}{S^{3/2}}.
\end{equation}
Since $b'(x)=a^{-1}(x)$, then for the risks difference (\ref{riskmain}) after some calculations we find 
 	\begin{equation*}
		\Delta (\theta) = c^2 -2c(d-1)\E_{\theta}\frac{1}{\sqrt{S}}
\end{equation*}
and we obtain a). The condition $\mathbf{C}$ implies
\begin{equation*}
		\Delta (\theta) = c^2-2c(d-1)a^*.
\end{equation*}
Obviously, the function $f(c)=c^2-2c(d-1)a^{*}$ is negative for any $c \in \left(0; 2(d-1)a^{*}\right)$, hence
$\Delta(\theta) \leq 0$ and we obtain b).
For c) we note that the extremum of function $f$ is achieved with $c=(d-1)a^{*}$. In this case we can write the explicit value for minimal gain :
$\Delta(\theta) \leq -(d-1)a^{*2}$. 

It remains to make sure that $\E_{\theta}(S^{-1/2})<\infty$. 
In case when the components of vector $X$ are independent, the finiteness of the expectation follows from \cite{Hudson1978}, since 
by Cauchy--Bunyakovskiy inequality $\E_{\theta}(S^{-1/2})\leq \bigl(\E_{\theta}(S^{-1})\bigr)^{1/2}$. For dependent components this result proved in \cite{Pchelintsev2013}.
Hence Theorem.

\bigskip
\section{Examples and corollaries}

In this section, we consider several examples of exponential distributions for which the  conditions for improved estimation hold.

\textbf{Example 1.} Let $X_1,...,X_d$ be independent Gaussian random variables with densities
	$$
 f_{\theta_i}(x)=\frac{1}{\sqrt{2\pi}}\exp\left\{-\frac{(x-\theta_i)^2}{2}\right\}\,, \quad i=1,...,d\,.
 $$
These functions belong to the class considered in this paper with 
	$$
 a(x)=1,\ \ b(x)=x,\ \ \psi(\theta_i)=-\frac{\theta_i^2}{2},\ \ k(x)=\frac{1}{\sqrt{2\pi}}\exp\{-x^2/2\}.
 $$
In this case $S=\norm{X}^2=\sum\limits_{i=1}^d X_i^2$, so the shrinkage estimates have the known form  \cite{Pchelintsev2013}
	\begin{equation}\label{delex1}
		\delta^*_{1}=X-\frac{c}{\|X\|}X\,.
	\end{equation}
To apply Theorem~\ref{Th1} one needs to check the condition $\mathbf{C}$.
It should be noted that we can represent any vector $X\sim\mathcal{N}_d(\theta,I)$ as $X=\theta + \xi$, where $\theta\in\Theta$ is an unknown mean vector, $\xi$ is a random vector from the standard normal distribution $\mathcal{N}_d(0,I)$, $I$ is $(d\times d)$ identity matrix.
 Scince the set $\Theta$ is bounded, then $\norm{\theta}\leq \tau$ for some positive constant $\tau$. Applying the Jensen inequality for expectation and the triangle inequality for norms, we obtain
 \begin{equation*}
 \E_\theta\frac{1}{\sqrt{S}}=\E_\theta\frac{1}{\|X\|} \geq \frac{1}{\E_\theta \|X\|}=\frac{1}{\E_\theta\|\theta + \xi\|} \geq \frac{1}{\|\theta\|+\E_\theta\|\xi\|}.    
 \end{equation*}
By Cauchy--Bunyakovskiy inequality $\E_\theta\|\xi\|\leq \sqrt{\E_\theta\|\xi\|^2}$. Taking into account that $\|\xi\|^2$ has the chi-squared distribution with $d$ degrees of
freedom and its expectation $\E_\theta\|\xi\|^2=d$, we get the inequality
\begin{equation*}
 \E_\theta\frac{1}{\sqrt{S}} \geq \frac{1}{\tau+\sqrt{d}}=:a^*    
 \end{equation*}
and, hence, the condition $\mathbf{C}$.

The risk of MLE
\begin{equation*}
	\cR(\theta,\delta^{0})=\E_\theta \| \theta-X\|^2=\E_\theta\|\xi\|^2=d\,.
\end{equation*}

Thus, we have the following result in Gaussian case.
\begin{cor}
    Let the observations vector $X \sim \mathcal{N}_d(\theta,I)$. Then for all $d \geq 2$ \\
   a) $\cR(\theta,\delta^*_1)=d+c^2-2c(d-1)\E_{\theta}(\norm{X}^{-1})$;\\
   b) the estimate \eqref{delex1} dominates the MLE and is minimax for all $0<c<2(d-1)a^{*}$, $a^{*}=(\tau+\sqrt{d})^{-1}$;\\
   c) the uniformly optimal choice of $c$ is 
   $c=(d-1)a^*$;\\
   d) in case c) $\cR(0,\delta^*)\to 0.5$ as $d\to\infty$.
\end{cor}

\begin{rem}
  The proof of d) is given in \cite{Pchelintsev2013}. The risk of MLE  $\cR(0,\delta^0)=d$ and tends to infinity as the dimension $d$ increases. The risk of James--Stein estimator  $\cR(0,\delta^{JS})= 2$. Thus, in the neighborhood of zero, the proposed estimate outperforms in accuracy to both the MLE and James--Stein estimator.
\end{rem}

Now we apply the estimate \eqref{delex1} for the multivariate normal distribution with correlated components. Let the vector $X \sim \mathcal{N}_d(\theta,V)$ with unknown mean vector $\theta$ and known or not covariance matrix $V$ such that its maximal eigenvalue 
\begin{equation}\label{lamda_max}
    \lambda_{\max}(V)\leq \lambda^{*}
\end{equation}
for some $\lambda^{*}<\infty$. 
The risk of MLE
$\cR(\theta,\delta^{0})=tr V$ (trace of matrix V) and we can check the condition $\mathbf{C}$ :
\begin{equation*}
 \E_\theta\frac{1}{\sqrt{S}} \geq \frac{1}{\E_\theta \|X\|} \geq \frac{1}{\|\theta\|+\E_\theta\|X-\theta\|}\geq \frac{1}{\tau+\sqrt{trV}}\geq \frac{1}{\tau+\sqrt{d\lambda^*}}=:a^*.    
\end{equation*}
Here we used Jensen, triangle and Cauchy--Bunyakovskiy inequalities, and estimate $trV\leq d\lambda^*$. Thus we have
\begin{cor}
    Let the observations vector $X \sim \mathcal{N}_d(\theta,V)$ and \eqref{lamda_max} hold. Then for all $d \geq 2$ \\
   a) $\cR(\theta,\delta^*_1)=tr V+c^2-2c(d-1)\E_{\theta}(\norm{X}^{-1})$;\\
   b) the estimate \eqref{delex1} dominates the MLE and is minimax for all $0<c<2(d-1)a^{*}$, $a^{*}=(\tau+\sqrt{d\lambda_{\max}})^{-1}$;\\
   c) the uniformly optimal choice of $c$ is 
   $c=(d-1)a^*$.
\end{cor}

\bigskip

\textbf{Example 2.} Let $X_1,...,X_d$ be independent random variables from the Gamma distribution with densities:
	$$
 f_{\theta_i}(x) = \frac{1}{\Gamma(\theta_i)}x^{\theta_i-1}e^{-x}\,, \quad x>0\,, \quad i=1,...,d\,.
 $$
Here $\Gamma(\cdot)$ is Gamma function. These densities belong to the class considered in this paper with 
	$$
 a(x)=x,\ \ b(x)=\ln x,\ \ \psi(\theta_i)=-\ln \Gamma(\theta_i),\ \ k(x)=x^{-1}\exp\{-x\}.
 $$
In this case $S=\norm{\ln X}^2=\sum\limits_{i=1}^d \ln^2 X_i$, so the shrinkage estimates have the form
 \begin{equation}\label{delex3}
		\delta_{2}^*=X-\frac{c}{\|\ln X\|}\ln X \,,\
	\end{equation}
 where	$\ln x = (\ln x_1, ... , \ln x_d).$
 
 Easy to see that the MLE for unknown vector of shape parameters $\theta$ is $\delta^0=X$ too. 
  Indeed, in general, the density function of the Gamma distribution has the form 
	\begin{equation*}
		f_{\mu,\lambda} (x) = \frac{\lambda^\mu}{\Gamma(\mu)}x^{\mu-1}e^{-\lambda x}\,.
	\end{equation*}
	The logarithmic likelihood function has the following form (since it consists of a single observation):
	\begin{equation*}
		\mathfrak{L}(\mu,\lambda)=\ln \left(f_{\mu,\lambda} (x)\right)=\mu \ln(\lambda)-\ln(\Gamma(\mu))+(\mu - 1)\ln x -\lambda x\,.
	\end{equation*}
	Solving the equation with respect to $\mu$
	\begin{equation*}
		\frac{\partial \mathfrak{L}}{\partial \lambda} = \frac{\mu}{\lambda}-x = 0\,, 
	\end{equation*}
	and using that $\lambda = 1$, $\mu=\theta_i$, $x=X_i$, we obtain desired view of the MLE. For quadratic risk  we have
	\begin{equation*}
		\cR(\theta,\delta^0)=\E_{\theta}\|X-\theta\|^2=\sum\limits_{i=1}^d \theta_i\,.
	\end{equation*}

 To apply Theorem~\ref{Th1} one needs to check the condition $\mathbf{C}$.
 		Let's estimate $\E_{\mu}(S)^{-1/2}$ from below by applying the Jensen and Cauchy--Buniakovskiy inequalities:
	\begin{equation}\label{S_gamma}
		\E_{\theta}\frac{1}{\sqrt{S}}=\E_{\theta}\left(\frac{1}{\|\ln X\|}\right)\geq \frac{1}{\sqrt{\sum\limits_{i=1}^d \E_{\theta} \left(\ln^2 X_i\right)}}.
	\end{equation}
Now we consider the expectation:
	\begin{equation}
		\E \bigl(\ln^2 x\bigr) = Var\left(\ln x\right)+\left(\E \left(\ln x\right)\right)^2.
	\end{equation}
	 In the calculations below, for simplicity, we will put $\theta_i=\mu$. We substitute $\ln\,x=y$, hence, $x=e^y$ and 
  $$
  f_\mu (y)=\frac{1}{\Gamma(\mu)}e^{(\mu-1)y - e^y}\,,
  $$
	$$\left(e^{(\mu-1)y - e^y}\right)'_\mu=y e^{(\mu-1)y-e^y}=y\Gamma (\mu) f_\mu (y).$$
Therefore,
	\begin{gather*}
		\E(y)=\int\limits_{\mathbb{R}} y f_\mu (y)\,dy=\int\limits_{\mathbb{R}}  \frac{y}{y\Gamma(\mu)}\left(e^{(\mu-1)y - e^y}\right)'_{\mu}\,dy=\frac{1}{\Gamma(\mu)}\int\limits_{\mathbb{R}}\left(e^{(\mu-1)y - e^y}\right)'_{\mu}\,dy= \\
		=\frac{1}{\Gamma(\mu)}\left(\int\limits_{\mathbb{R}} e^{(\mu-1)y - e^y}\,dy\right)'=\frac{\Gamma'(\mu)}{\Gamma(\mu)}=\psi(\mu)\,,\
	\end{gather*}
where $\psi(\mu)$ is a digamma function.
	Further, using the same changing of variables, we find
	\begin{gather*}
		Var(y)=\E\Bigl(y-\E y\Bigr)^2=\E\Bigl(y-\psi (\mu)\Bigr)^2=\int\limits_\mathbb{R} \left(y-\psi (\mu)\right)^2 f_\mu (y)\,dy\\
		=\int\limits_\mathbb{R} y^2 f_\mu (y)\,dy -\psi^2(\mu)=\int\limits_\mathbb{R} \frac{y }{\Gamma(\mu)}\left(e^{(\mu-1)y - e^y}\right)'_{\mu}\,dy-\psi^2(\mu) \\
		=\frac{1}{\Gamma(\mu)}\left(\int\limits_\mathbb{R} \left( e^{(\mu-1)y - e^y}\right)'_{\mu}\,dy\right)'_{\mu}-\psi^2(\mu)=\frac{1}{\Gamma(\mu)}\left(\int\limits_\mathbb{R}  e^{(\mu-1)y - e^y}\,dy\right)''_{\mu}-\psi^2(\mu) \\
		=\frac{\Gamma''(\mu)}{\Gamma(\mu)}-\left(\frac{\Gamma'(\mu)}{\Gamma(\mu)}\right)^2=\frac{\Gamma(\mu)\Gamma'(\mu)-(\Gamma'(\mu))^2}{\Gamma^2(\mu)}=\left(\frac{\Gamma'(\mu)}{\Gamma(\mu)}\right)'=\psi'(\mu)=\psi^{(1)}(\mu)\,,
	\end{gather*}
where $\psi^{(1)}(\mu)$ is the trigamma function.
	
	Thus, $\E_{\mu} (\ln^2 x)=\psi^{(1)}(\mu)+(\psi(\mu))^2$, and we need to estimate from above this sum.	
	An estimate for $\psi(\cdot)$ is proposed in \cite{8}:
	$$\frac{1}{2\mu}<\ln \mu -\psi(\mu)<\frac{1}{\mu}\,, \quad \mu>0.$$
	To prove the inequality, the function $H(x)=x(\ln\,x -\psi(x))$ is defined. It is proved in [9] that this function is strictly decreasing and strictly convex on $(0,\infty)$. 
	Let us also use 
	$$
 H(x)=\frac{1}{2}+\frac{1}{12x}+\frac{\nu}{120x^3}.
 $$
	The proof is given for the case of the lower bound, because it is necessary for further reasoning.
	We define integral representations for $\ln x$, $\psi(x)$, and $\frac{1}{2x}$
	\begin{equation*}
		\ln\,x -\psi(x)-\frac{1}{2x}=\int\limits_0^\infty e^{-xt}\left(\frac{e^t}{e^t -1}-\frac{1}{t}-\frac{1}{2}\right)\,dt.
	\end{equation*}
	Function
	$$
 H(t)=2t(e^t-1)\left(\frac{e^t}{e^t -1}-\frac{1}{t}-\frac{1}{2}\right)=t e^t-2 e^t + t+2
 $$
	is positive for $t>0$. Moreover, $H'(t)=t e^t -e^t+1$, $H'(0)=0$ and $H'(t)=t e^t>0$. $H(t)$ increases by $(0,\infty)$. Using this fact and $H(0)=0$, we conclude that $H(t)>0$ and obtain the lower bound.

Now let's estimate $\psi^{(1)}(\mu)$, using Stirling expansion
	$$
	\ln \Gamma (\mu) = \left(\mu - \frac{1}{2}\right)\ln(\mu) -\mu +\frac{1}{2}\ln(2\pi)+\sum\limits_{n=1}^m \frac{B_{2n}}{2n(2n-1)\mu^{2n-1}} + O(\mu^{-2m-1})\,,\
	$$
	where $B_{2n}$ are Bernoulli numbers,
	\begin{gather*}
		\psi^{(0)}(\mu)=\frac{d}{d\mu} \ln \Gamma (\mu) = \ln \mu - \frac{1}{2\mu} - \sum\limits_{n=1}^m \frac{B_{2n} \mu^{-2n}}{2n} + O(\mu^{-2m-1})\,, \\
		\psi^{(1)}(\mu) = \frac{d}{d\mu} \psi^{(0)}(\mu)= \frac{1}{\mu} + \frac{1}{2\mu^2} + \sum\limits_{n=1}^m \frac{B_{2n}}{\mu^{2n+1}} + O(\mu^{-2m-1}).
	\end{gather*}
	From the last relation, the inequality directly follows 
	$$
	\psi^{(1)}(\mu)<\frac{1}{\mu}+\frac{1}{2 \mu^2}+\frac{1}{6 \mu^3}.
	$$
	Using inequalities for the polygamma functions to estimate the expectation, we get
	\begin{gather*}
		\E (\ln^2 x)=\psi^{(1)}(\mu) + (\psi(\mu))^2<\frac{1}{\mu}+\frac{1}{2\mu^2}+\frac{1}{6\mu^3}+\left(\ln(\mu)-\frac{1}{2\mu}\right)^2 \\
		=\frac{12\mu^2-12\mu^2 \ln(\mu)+9\mu+12\mu^3\ln^2(\mu)+2}{12\mu^3} =: M(\mu).
	\end{gather*}

\hspace{1mm}
\begin{center}
\includegraphics[width=10cm]{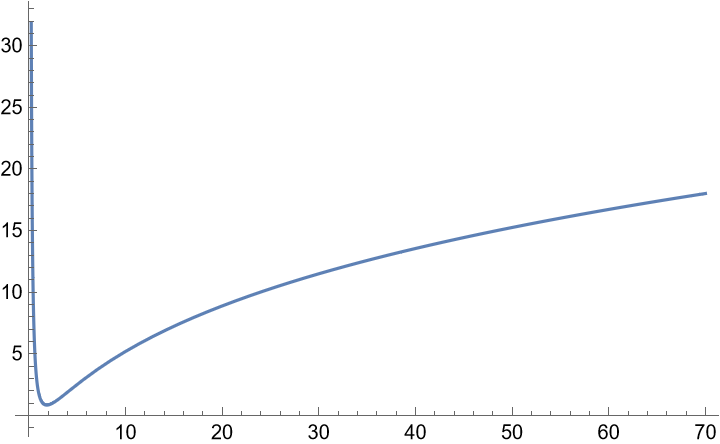}\\
{\footnotesize Figure 1. Graphic of function $M(\mu)$} \end{center}\label{fig1} \hspace{3mm}
 

 As we can see in Fig.~\ref{fig1}, function $M(\mu)$ has exclusively a global minimum at $\mu=1.865$. Accordingly, two cases arise.
	
	If $\mu \leq 1.865$, we define $\mu^{*}$ to estimate $M(\mu)$ from above: $ \Theta$ is bounded set, then there exists $\mu^{*}:= \inf\{\norm{\theta}\,:\,\theta \in \Theta\}$, hence, any $\mu \geq \mu^{*}$ and
$$
\E(\ln^2 x)<\frac{12\mu^{*2}-12\mu^{*2} \ln(\mu^{*})+9\mu^{*}+12\mu^{*3}\ln^2(\mu^{*})+2}{12\mu^{*3}}=:m_1\,.
$$
	
	If $\mu > 1.865$, we define $\mu^{**}$ to estimate $M(\mu)$ from above:  $\mu^{**}:= \sup\{\norm{\theta}\,:\,\theta \in \Theta\}$, hence any $\mu \leq \mu^{**}$ and
\begin{equation*}
		\E(\ln^2 x)<\frac{12\mu^{**2}-12\mu^{**2} \ln(\mu^{**})+9\mu^{**}+12\mu^{**3}\ln^2(\mu^{**})+2}{12\mu^{**3}}=:m_2.
\end{equation*}
	Then, using this estimates in \eqref{S_gamma} we obtain the condition $\mathbf{C}$ with 
	$a^{*}=(d m)^{-1/2}$, where $m=\max(m_1, m_2)$.
		So, we proved 
\begin{cor}
    Let the observations vector $X \sim Gamma(\theta,\mathbf{1})$. Then 
    for all $d \geq 2$ \\
   a) $\cR(\theta,\delta^*_2)=\cR(\theta,\delta^0)+c^2-2c(d-1)\E_{\theta}(\norm{\ln X}^{-1})$;\\
   b) the estimate \eqref{delex3} dominates the MLE and is minimax for all $0<c<2(d-1)a^{*}$, $a^{*}=(md)^{-1/2}$;\\
   c) the uniformly optimal choice of $c$ is 
   $c=(d-1)a^*$.
\end{cor}

\bigskip

\textbf{Example 3.} Let $X_1,...,X_d$ be independent random variables from exponential distribution with densities:
	$$f_{\theta_{i}}(x) = \frac{1}{\theta_{i}} \cdot e^{-\frac{x}{\theta_{i}}}\,, \quad x>0\,, \quad i=1,...d\,.$$
Putting $\lambda_i=1/\theta_i$ and $y=-x$, we can represent density $f_{\lambda_{i}}(y)$ in the form \eqref{den} if
$$
a(y)=1,\ \ b(y)=y,\ \ \psi(\lambda_i)=-\ln \lambda_i,\ \ k(y)=1\,. 
$$
But in general for this distribution the identity \eqref{lemStein} does not true. However proposed method can be apply for this case. Direct calculations give that the MLE for unknown mean vector $\theta=(\theta_1,...,\theta_d)$ is $\delta^0=-Y=X$, and its quadratic risk  
	\begin{equation*}
		\cR(\theta,\delta^0)=\E_{\theta}\|X-\theta\|^2=\E_{\theta} \sum\limits_{i=1}^d (X_i-\theta_i)^2 = \sum\limits_{i=1}^d \theta_i^2=\norm{\theta}^2\,.
	\end{equation*}
Since $S=\sum\limits_{i=1}^d X_i^2$, then the shrinkage estimate for unknown mean vector $\theta$ be defined as 
	\begin{equation}\label{delta3}
		\delta_{3}=X - \frac{c}{\|X\|}X\,.
	\end{equation}
The quadratic risk  of this estimate 
	\begin{equation*}
		\cR(\theta,\delta^*_3)=\cR(\theta,\delta^0)+c^2+2\sum\limits_{i=1}^d \E_{\theta} (X_i-\theta_i)g_i(X) \,,
	\end{equation*}
where $g_i(x)=-cx_i/\norm{x}$. Integrating by parts, for expectation in this sum we have
$$
\E_{\theta} (X_i-\theta_i)g_i(X)=\int...\int (x_i-\theta_i)g_i(x)\prod_{j=1}^d f_{\theta_j}(x_j)dx_j=\theta_i\E_{\theta} [X_ig'_i(X)]\,.
$$
This equality and \eqref{deriv_g} imply that the risks difference  
	\begin{equation*}\label{risk3}
		\Delta (\theta)\leq c^2 - 2c(d-1)\theta^*\E_{\theta}\frac{1}{\|X\|} \,,
	\end{equation*}
 where $\theta^*:=\inf\{\norm{\theta}\,:\,\theta \in \Theta\}$. 
To estimate last expectation from below we apply the Jensen and Cauchy--Bunyakovskiy inequalities:
	\begin{equation*}\label{ineq2}
		\E_{\theta}\frac{1}{\|X\|}=\frac{1}{\E_{\theta}\norm{X}} \geq \frac{1}{\sqrt{\sum\limits_{i=1}^d \E_{\theta} X_i^2}}\,.
	\end{equation*}
It easy to see that $\E_{\theta} X_i^2=VarX_i+(\E_{\theta} X_i)^2=2\theta_i^2$. 
	Let $\theta^{**}:=\sup\{\norm{\theta}\,:\,\theta\in \Theta\}$.
	 Then we find that
	\begin{equation*}
		\Delta (\theta) \leq c^2-2c(d-1)a^*\,,
	\end{equation*}
	where $a^{*}=\theta^{*}/(\sqrt{2}\theta^{**})$.
Hence, we have
\begin{cor}
    Let the observations vector $X \sim Exp(\theta)$. Then 
   for all $d \geq 2$ \\
   a) $\cR(\theta,\delta^*_3)=\norm{\theta}^2+c^2-2c\sum_{i=1}^d\theta_i\E_{\theta}[X_ig'_i(X)]$;\\
   b) the estimate \eqref{delex3} dominates the MLE and is minimax for all $0<c<2(d-1)a^{*}$, $a^{*}=\theta^{*}/(\sqrt{2}\theta^{**})$;\\
   c) the uniformly optimal choice of $c$ is 
   $c=(d-1)a^*$.
\end{cor}

\bigskip
\section{Numerical simulation}	
In this section we give  a numerical simulation to compare the empirical quadratic risks of three estimates: MLE $\delta^0$, the James-Stein $\delta^{JS}$ or Hudson $\delta^{H}$ estimate, and the proposed improved estimates $\delta^{*}$ in Sect. 2-3. We will simalate the estimate $\delta^{*}$ with optimal $c=(d-1)a^*$.
	
We generate observations vector $X$ of size $d$ from normal, Gamma and exponential distribution. Parameters $\theta_i$ varying according to different formulas.  The empirical risk of any estimate $\delta$ we define as
\begin{equation}
    \hat{R}(\theta,\delta)= \frac{1}{1000} \sum\limits_{l=1}^{1000} \|\delta_l -\theta\|^2\,,
\end{equation}
where $\delta_l$ is $l$-th replication of estimate $\delta$.
To compare the empirical risks we calculate the ratios of them, i.e.
$$
R^{JS}=\frac{\hat{R}(\theta,\delta^0)}{\hat{R}(\theta,\delta^{JS})}\,,\quad R^{H}=\frac{\hat{R}(\theta,\delta^0)}{\hat{R}(\theta,\delta^{H})}\,,\quad
R^{*}=\frac{\hat{R}(\theta,\delta^0)}{\hat{R}(\theta,\delta^{*})}\,.
$$

\textbf{Normal distribution.} For the normal distribution, we consider three cases: in the case A),  the mean vector $\theta$ is not large; in the case B), the observations are i.i.d normal (0,1); in the case C),  the mean vector $\theta$ is  large. Here we apply Corollarry 3.1 and Remark 3.2.
The results of the simulation are presented in Table 1.
\begin{table}[H]
		\centering
		\begin{tabular}{|c|c|c|c|c|c|}
			\hline
			$d$&2&10&50&100&500  \\
   \hline \hline
\multicolumn{6}{|l|}{A) $X_i\sim\mathcal{N}(\theta_i, 1)$, $\theta_i=(i+1)/d$, $i=1,...,d$} \\
			\hline
			$R^{JS}$& 1.00 & 2.05 & 3.41 & 3.64 & 3.93\\
			\hline
			$R^{*}$& 1.18 & 2.25 & 3.01 & 3.13 & 3.26\\
			\hline
\hline
\multicolumn{6}{|l|}{B) $X_i\sim\mathcal{N}(\theta_i, 1)$, $\theta_i=0$, $i=1,...,d$} \\
			\hline
			$R^{JS}$& 1.00 & 4.97 & 24.64 & 49.36 & 242.83\\
			\hline
			$R^{*}$& 2.78 & 18.26 & 97.84 & 193.35 & 974.16\\
			\hline
   \hline
\multicolumn{6}{|l|}{C) $X_i\sim\mathcal{N}(\theta_i, 1)$, $\theta_i=5+id$, $i=1,...,d$} \\
			\hline
			$R^{JS}$& 1.00 & 1.0002 & 1.00 & 1.00 & 1.00\\
			\hline
			$R^{*}$& 1.0058 & 1.0003 & 1.00 & 1.00 & 1.00\\
			\hline
		\end{tabular}
		\caption{Empirical risks ratios of estimates for normal distribution case}
		\label{tab:1}
	\end{table}

 Note that the proposed estimate dominates the MLE for all dimensions $d\geq 2$. In the case A) the mean square accuracy of shrinkage estimates more then for MLE in 2-3 times. For not large $d$ the proposed estimate also outperforms the James--Stein estimate. 
 In the case B) we have superdomination effect for the proposed estimate over MLE and James--Stein estimate. Case C) illustrates the fact that James--Stein-type estimators have an advantage for small parameters. 
 It can be seen that for large parameters $\theta$ the gain in accuracy becomes insignificant. In such cases the influence of the bias is reduced and the estimation becomes close to MLE.

 \textbf{Gamma distribution.}
 Now we simulate the observations from the Gamma distribution. The Hudson modification applicable to the Gamma distribution is used instead of the James--Stein estimate. We consider two cases: in the case D),  the mean vector $\theta$ is not large and in the case E), the mean vector $\theta$ is  large. Here we apply Corollarry 3.4.
The results of the simulation are presented in Table 2.
\begin{table}[H]
		\centering
		\begin{tabular}{|c|c|c|c|c|c|}
			\hline
			$d$&2&10&50&100&500  \\
   \hline \hline
\multicolumn{6}{|l|}{D) $X_i\sim Gamma(\theta_i, 1)$, $\theta_i=(i+1)/d$, $i=1,...,d$} \\
			\hline
			$R^{H}$& 1.00 & 1.21 & 1.12 & 1.11 & 1.09\\
			\hline
			$R^{*}$& 1.33 & 1.12 & 1.06 & 1.05 & 1.04\\
			\hline
   \hline
\multicolumn{6}{|l|}{E) $X_i\sim Gamma(\theta_i, 1)$, $\theta_i=5+id$, $i=1,...,d$} \\
			\hline
			$R^{H}$& 1.00 & 1.0007 & 1.00 & 1.00 & 1.00\\
			\hline
			$R^{*}$& 1.009 & 1.0008 & 1.00 & 1.00 & 1.00\\
			\hline
		\end{tabular}
		\caption{Empirical risks ratios of estimates for Gamma distribution case}
		\label{tab:2}
	\end{table}
     Here it is also possible to conclude that the proposed estimate enables greater accuracy for not large parameter values, and the gain in the accuracy of shrinkage estimates turns out to be small.

\bigskip

\textbf{Exponential distribution.} 
We simulate the observations from the exponential distribution. For this distribution the James--Stein or Hudson estimate cannot construct. We also consider two cases: in the case F),  the mean vector $\theta$ is not large and in the case G), the mean vector $\theta$ is  large. Here we apply Corollarry 3.5.
The results of the simulation are presented in Table 3.
\begin{table}[H]
		\centering
		\begin{tabular}{|c|c|c|c|c|c|}
			\hline
			$d$&2&10&50&100&500  \\
   \hline \hline
\multicolumn{6}{|l|}{F) $X_i\sim Exp(\theta_i)$, $\theta_i=(i+1)/d$, $i=1,...,d$} \\
			\hline
			$R^{*}$& 5.03 & 2.17 & 1.38 & 1.36 & 1.32\\
			\hline
   \hline
\multicolumn{6}{|l|}{G) $X_i\sim Exp(\theta_i)$, $\theta_i=5+id$, $i=1,...,d$} \\
			\hline
			$R^{*}$& 2.47 & 1.04 & 1.01 & 1.36 & 1.31\\
			\hline
		\end{tabular}
		\caption{Empirical risks ratios of estimates for exponential distribution case}
		\label{tab:3}
	\end{table}

As we can see for this case the proposed estimate dominates the MLE for all dimensions $d\geq 2$. The mean square accuracy of shrinkage estimate more then for MLE in 2-5 times if parameter values are not large, and in 1.3-1.5 times else.

\bigskip
\section{Conclusion}
This paper proposes an improved estimation procedure for exponential distribution parameters. The result about improvement based on the assumption that the vector of unknown parameters belong to a compact set. 

Due to numerical simulations, the proposed biased estimate actually has smaller risk than the maximum likelihood estimate.

\bigskip
\section{Acknowledgment}
This work has been supported by RSF, project no 20-61-47043.

\end{document}